\documentclass[11pt]{article}
\usepackage[utf8]{inputenc}
\usepackage[letterpaper,margin=1in]{geometry}
\usepackage{amstext,amsfonts,amssymb,amscd,amsbsy,amsmath,tikz,mathrsfs,amsthm,float,graphicx}
\usepackage{caption,subcaption}
\usepackage{mathtools}

\newtheorem{theorem}{Theorem}[section]

\newtheorem{corollary}[theorem]{Corollary}
\newtheorem{lemma}[theorem]{Lemma}
\newtheorem{proposition}[theorem]{Proposition}

\newtheorem{question}[theorem]{Question}

\newcommand{\Z}{\mathbb{Z}}

\newcommand{\R}{\mathbb{R}}

\DeclareMathOperator{\sat}{sat}
\DeclareMathOperator{\ex}{ex}
\DeclareMathOperator{\KR}{KR}
\DeclareMathOperator{\argmin}{argmin}

\providecommand{\keywords}[1]
{
  \small	
  \textbf{\textit{Keywords---}} #1
}

\title{Generalized saturation problems for cliques, paths, and stars}
\author{Jamie Radcliffe$^1$ \and Adam Volk$^1$}
\date{%
    $^1$University of Nebraska-Lincoln\\[2ex]%
    \today
}

\begin{document}

\maketitle

\begin{abstract}
    A graph $G$ is $F$-\textit{saturated} if it does not contain any copy of $F$, but the addition of any missing edge in $G$ creates at least one copy of $F$. Inspired by work of Alon and Shikhelman regarding a similar question for $F$-free graphs, Kritschgau, Methuku, Tait, and Timmons introduced the parameter of $\sat_H(n,F)$ to denote the minimum number of copies of some subgraph $H$ in an $F$-saturated graph on $n$ vertices. In this paper, we address this generalized saturation problem with special focus on $\sat_{K_r}(n,S_t)$ and $\sat_{S_r}(n,S_t)$. This relates to recent work by Chakraborti and Loh regarding $\sat_{K_r}(n,K_t)$ and by Ergemlidze, Methuku, Tait, and Timmons regarding $\sat_{S_r}(n,K_t)$. We also provide some results regarding paths and arbitrary trees. 
\end{abstract}


\section{Introduction}
A major object of focus in extremal graph theory is the extremal number $\ex(n,F)$, which denotes the maximum number of edges in an $F$-free graph on $n$ vertices. The study of such values dates back to Tur\'an in 1941 \cite{T}. See \cite{S} for a survey by Sidorenko with some results in the field. More recently, Alon and Shikhelman \cite{AS} introduced the generalized extremal number $\ex_H(n,F)$ which gives the maximum number of copies of $H$ among $F$-free graphs on $n$ vertices. If we take  $H$ to be a single edge, we recover the original extremal number. 

In addition to maximizing the number of copies of a given subgraph among $F$-free graphs, it is also natural to try to understand the minimum. To avoid trivialities, we say that a graph $G$ is $F$-\textit{saturated} if $G$ does not contain any copy of $F$ as a subgraph, but the addition of any missing edge creates at least one copy of $F$. We do not consider induced subgraphs in this setting although it is possible to do so. (For instance, see \cite{MS} for a consideration of induced saturation.) The minimum number of edges in an $F$-saturated graph on $n$ vertices is denoted $\sat(n,F)$. The case where $F$ is the complete graph was solved by Erd\H{o}s, Hajnal, and Moon \cite{EHM} and is stated here.

\begin{theorem}[Erd\H{o}s, Hajnal, and Moon, 1964]\label{theorem_ehm}
For every $n\geq t\geq 2$, 
$$\sat(n,K_t)=(n-t+2)(t-2)+\binom{t-2}{2}.$$ The graph $K_{t-2}+\overline{K}_{n-t+2}$ is the unique extremal example.
\end{theorem}

Saturation numbers have been well-studied and a collection of some of these results can be found in \cite{FFS}. Inspired by the generalization of Alon and Shikhelman in the $F$-free setting, Kritschgau, Methuku, Tait, and Timmons \cite{KMTT} introduced the generalized saturation number $\sat_H(n,F)$ to denote the minimum number of copies of $H$ in an $F$-saturated graph on $n$ vertices. In addition to proving general results, they focused on the cases where at least one of $H$ and $F$ was a clique or cycle. Extending a result of Kritschgau et al. and proving a conjecture from that paper, Chakraborti and Loh \cite{CL} showed the following. 

\begin{theorem}[Chakraborti and Loh, 2019]
    For every $t> r \geq 2$, there exists a constant $n_{r,t}$ such that, for all $n\geq n_{r,t}$, we have 
    $$\sat_{K_r}(n,K_t)=(n-t+2){{t-2}\choose {r-1}}+{{t-2}\choose {r}}.$$
    Furthermore, for $n$ sufficiently large, the (complete) split graph $K_{t-2}+\overline{K}_{n-t+2}$ is the unique extremal example. 
\end{theorem}

Chakraborti and Loh asked if the split graph minimizes the number of copies of any $F$ among $K_t$-saturated graphs. Using stars as their choice of $F$,
Ergemlidze, Methuku, Tait, and Timmons \cite{EMTT} proved that this is not the case. Furthermore, the split graph was far from optimal. 

In this paper, we look at the other variations involving stars and cliques, namely $\sat_{K_r}(n,S_t)$ and $\sat_{S_r}(n,S_t)$. We also consider $\sat_{S_r}(n,K_t)$ in a more restricted setting as well as generalized saturation numbers involving paths.

\subsection{Notation and organization}
In general our notation follows that of Bollob\'{a}s \cite{BB}. In particular, we write $G_1 +G_2$ to denote the join of two graphs $G_1$ and $G_2$. That is, we take the disjoint union of the two graphs and add the edge $uv$ for every vertex $u$ in $G_1$ and $v$ in $G_2$. We similarly write $G_1 \cup G_2$ for the disjoint union of $G_1$ and $G_2$ with no added edges. For a given vertex $v$, we write $N(v)$ to denote the neighborhood of $v$ and $N[v]$ to denote the closed neighborhood of $v$. That is, $N[v] = N(v) \cup\{v\}$. We write $(n)_k$ to denote the falling factorial $n(n-1)\cdots(n-k+1)$ and will make use of the generalized binomial coefficient which is defined to be ${{x}\choose{k}} = \frac{1}{k!}x(x-1)\cdots(x-k+1)$ for an arbitrary real number $x$ and integer $k$.

Given a graph $G$, we let $k_r(G)$ and $s_r(G)$ denote the number of copies of cliques $K_r$ and stars $S_r$ in $G$ respectively. Using the notation of Gallian \cite{JG}, we let $S_r$ denote a star on $r+1$ vertices and $r$ edges. 

The paper is organized as follows. In Section \ref{sec_clique_star} and Section \ref{sec_star_in_star} we count cliques and stars respectively in star-saturated graphs. In Section \ref{sec_star_clique}, we briefly discuss the problem of counting stars in clique-saturated graphs. In particular, we comment on a recent result of Ergemlidze, Methuku, Tait, and Timmons \cite{EMTT} and state a result regarding clique-saturated graphs with linear maximum degree. In Section \ref{sec_paths}, we provide an asymptotic result regarding the minimum number of paths in clique-saturated graphs. We also use a result of Kaszonyi and Tuza \cite{KT} to show that for $n$ sufficiently large, the minimum number of copies of $K_r$ in path-saturated graphs is in fact $0$ for all $r\geq 3$. We finish in Section \ref{sec_general} with a brief discussion of some more general results regarding cliques and trees.

\subsection{Graph constructions}\label{Constructions}
In the sections that follow, we will make use of certain graphs that are regular or almost regular. The first construction concerns regular multipartite graphs. Before presenting the construction, we state a theorem of Hoffman and Rodger \cite{HoffRog}. 

\begin{theorem}[Hoffman and Rodger, 1992]\label{HoffRog}
    Given a complete multipartite graph $K$, $\chi'(K) = \Delta(K)$ if and only if it is not overfull. Here $\chi'$ denotes the chromatic index of $K$, and we say that a graph $G$ is \textit{overfull} if $|E(G)| > \Delta(G)\left\lfloor\frac{|V(G)|}{2}\right\rfloor$. In particular, the complete $r$-partite graph $K_{a,\dots,a}$ is overfull if and only if $ar$ is odd.
\end{theorem}
\begin{proposition}\label{hoff_const}
    Let $a,r$ be positive integers. If $ar$ is even, there exists a $k$-regular spanning subgraph of the $r$-partite graph $K_{a,\dots,a}$ for all $k \leq a(r-1)$. 
    \begin{proof}
        Consider the complete $r$-partite graph $K_{a,\dots,a}$. If $ar$ is even, then $K_{a,\dots,a}$ is not overfull. By Theorem \ref{HoffRog}, the chromatic index of this graph is equal to the maximum degree. That is, we can give a proper edge coloring using $a(r-1)$ colors. In particular, the color classes form a $1$-factorization of $K_{a,\dots,a}$. We delete perfect matchings until we are left with a $k$-regular subgraph for any $k\leq a(r-1)$. 
    \end{proof}
\end{proposition}

For our next construction, given $a<b$, we let $R_{a,b}$ denote a graph on $b$ vertices that is as close to $a$-regular as possible. More specifically, when $ab$ is even, $R_{a,b}$ is $a$-regular. When $ab$ is odd, $R_{a,b}$ has one vertex of degree $a-1$ and $b-1$ vertices of degree $a$.

\begin{lemma}\label{Rab_lemma}
    An $R_{a,b}$ exists if and only if $b\geq a+1$.
    \begin{proof}
        It is well known that $a$-regular graphs exist on $b$ vertices if and only if $b\geq a+1$ and $ab$ is even. When $b\geq a+1$ and $ab$ is odd, we can obtain a graph on $b$ vertices that is $a$-regular with the exception of one vertex of degree $a-1$ in the following manner. Label the vertices $0,1,\dots,b-1$. Add edges between vertices with labels $i$ and $j$ if and only if $$|i-j|\leq \frac{a-1}{2}\mod{b}.$$ Finally for $1 \leq i \leq \frac{b-1}{2}$, add an edge between the vertices labeled $i$ and $i+\frac{b-1}{2}\mod{b}$. We aren't reusing any edges since $b\geq a +1$, and the degree of every vertex increases by $1$, except for the vertex labeled $0$. Thus every vertex has degree $a$ with the exception of one vertex of degree $a-1$, and $R_{a,b}$ exists. Necessity is clear because a graph can't have any vertex of degree $a$ if there are fewer than $a+1$ vertices.
    \end{proof}
\end{lemma}

Utilizing the regularity of $R_{a,b}$, we define some candidate extremal graphs for $\sat_{S_r}(n,S_t)$ that will be studied in Section \ref{sec_star_in_star}. They are essentially the disjoint union of a clique and a regular graph. To be precise, for $m\leq t-1$ we let 
\[
    \KR_{t,n}(m) = \begin{cases}
                        K_m \cup R_{t-1,n-m} & \text{if $(t-1)(n-m)$ is even} \\
                        K_m \cup R_{t-1,n-m} + e & \text{if $(t-1)(n-m)$ is odd},
                    \end{cases}
\]
where in the second case $e$ is an edge between the vertex of degree $t-2$ in $R_{t-1,n-m}$ and an arbitrary vertex of the clique $K_m$.

\section{Cliques in star-saturated graphs}\label{sec_clique_star}
We begin this section with a structural lemma regarding star-saturated graphs, which was observed by Kaszonyi and Tuza \cite{KT}. We include the proof for completeness.
\begin{lemma}\label{lemma_2.1}
    Let $G$ be an $S_t$-saturated graph on $n$ vertices. Then the maximum degree of $G$ is $t-1$, and all vertices of degree less than $t-1$ form a clique.
    \begin{proof}
        Since $G$ is $S_t$-saturated, it must be $S_t$-free. Thus the maximum degree of $G$ is at most $t-1$. Adding any missing edge can increase the maximum degree by at most $1$. Since $G$ is $S_t$-saturated, the new edge must force the maximum degree to become $t$. Thus $G$ has maximum degree exactly $t-1$. 
        
        To prove the second part of the statement, suppose that $u$ and $v$ are vertices in $G$ with degrees $d(u) < t- 1$ and $d(v)< t-1$. If $u$ is not adjacent to $v$, then the addition of the edge $uv$ will only increase their degrees, and the resulting graph will still have maximum degree less than $t$, a contradiction to $G$ being $S_t$-saturated. Therefore all vertices of degree less than $t-1$ must form a clique in $G$. 
    \end{proof}
\end{lemma}

The following result allows us to completely determine the minimum number of copies of $K_r$ in an $S_t$ saturated graph for all $r\geq 3$ and $t\geq 3$ when $n$ is sufficiently large. Note that for $t =1$ we need $n\geq 2$, and for $t =2$ we need $n \geq 3$ for $S_t$-saturated graphs to exist. In particular, there is a unique $S_2$-saturated graph on $n$ vertices for all $n\geq 3$. When $n$ is even, the graph is a collection of disjoint edges. When $n$ is odd, it is a collection of disjoint edges and an isolated vertex. $S_1$-saturated graphs are simply independent sets on at least $2$ vertices.

\begin{proposition}\label{prop_K3_iff}
    Let $t\geq 3$. There exists an $S_t$-saturated graph on $n$ vertices that is $K_3$-free if and only $n\geq 2t-2$. 
    \begin{proof}
        If $n$ is even and $n\geq 2t-2$, there exists an $(t-1)$-regular bipartite graph on $n$ vertices. This is certainly $S_t$-saturated. When $n$ is odd, consider $G\cup K_1$ where $G$ is a $(t-1)$-regular bipartite graph on $n-1$ vertices. The resulting graph is still $S_t$-saturated and bipartite, hence $K_3$-free.
        
        For the other direction, let $G$ be $S_t$-saturated and $K_3$-free. If $G$ has no vertices of degree $t-1$, then by Lemma \ref{lemma_2.1} $G$ is complete and contains triangles. This means there is a vertex $v \in G$ with degree $t-1$. For any two vertices $x,y\in N(v)$, at least one of them must have degree $t-1$. Otherwise $x$ is adjacent to $y$ and $v,x,y$ form a copy of $K_3$. Without loss of generality, the degree of $x$ is $t-1$. Since $G$ is $K_3$-free, $N(v)\cap N(x) = \varnothing$, and $N(x)$ contains $t-2$ vertices outside of $N[v]$. Therefore $n\geq 2t-2$.
    \end{proof}
\end{proposition}


The following is immediate.
\begin{corollary}
    For all $r\geq 3$ and all $t\geq 3$, $$\sat_{K_r}(n,S_t) = 0$$ for $n\geq 2t-2$. When $t =1,2$, $$sat_{K_r}(n,S_t) = 0$$ when $n\geq 2,3$ respectively.
\end{corollary}

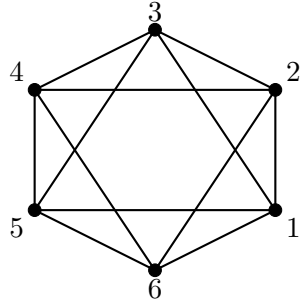
\begin{figure}[H]
    \centering
    \begin{tikzpicture}[scale =0.8]
        \draw[fill=black] (2,-1) circle (3pt);
        \draw[fill=black] (2,1) circle (3pt);
        \draw[fill=black] (0,2) circle (3pt);
        \draw[fill=black] (-2,1) circle (3pt);
        \draw[fill=black] (-2,-1) circle (3pt);
        \draw[fill=black] (0,-2) circle (3pt);
        \node at (2.3,-1.3) {1};
        \node at (2.3,1.3) {2};
        \node at (0,2.3) {3};
        \node at (-2.3,1.3) {4};
        \node at (-2.3,-1.3) {5};
        \node at (0,-2.3) {6};
        \draw[thick] (2,-1) -- (2,1) -- (0,2) -- (-2,1) -- (-2,-1) -- (0,-2) -- (2,-1);
        \draw[thick] (2,-1) -- (0,2) -- (-2,-1) -- (2,-1);
        \draw[thick] (2,1) -- (-2,1) -- (0,-2) -- (2,1);
    \end{tikzpicture}
    \caption{An example of a $K_4$-free, $S_5$-saturated graph on fewer than $8$ vertices}
    \label{fig:has_k4} 
\end{figure}

Proposition \ref{prop_K3_iff} gives us a cutoff for the values of $n$ which require $S_t$-saturated graphs to contain at least one copy of $K_3$. It tells us, for instance, that any $S_5$-saturated graph that is $K_3$-free must have at least $8$ vertices. Although there does not exist such a graph on $6$ vertices, the graph in Figure \ref{fig:has_k4} is an example of an $S_5$-saturated graph on $6$ vertices that, although not $K_3$-free, is $K_4$-free.

These results, together with this observation, lead us to the following question.
\begin{question}\label{ques2}
    Given a fixed $r \geq 3$ and $t\geq 3$, for which values of $n$ do there exist $S_t$-saturated graphs that are $K_{r+1}$-free? Furthermore, what are the possible values of the number of copies of $K_r$ in $S_t$-saturated graphs on $n$ vertices?
\end{question}

We provide a partial answer to the first question, focusing especially on the existence of $S_t$-saturated graphs that are $r$-partite. Due to the triviality of the cases where $t =1,2$, our general results consider $t\geq 3$.

\begin{theorem}\label{theorem_partite}
    Let $r\geq 3$ and $t\geq 3$ be fixed. There exists an $n$-vertex, $r$-partite, $S_t$-saturated graph if 
    \begin{equation}\label{ineq3}
        n\geq \max \left( t+1,\min_{0\leq c \leq r-2}\left\{(r-c)\left\lceil\frac{t-1}{r-c-1}\right\rceil + r-c\right\}\right).
    \end{equation}  
    \begin{proof}
        Suppose (\ref{ineq3}) holds for some $0\leq c\leq r-2$. Let $a,b$ be non-negative integers such that $n = a(r-c) + b$ with $b < r-c$. That is, $a = \lfloor \frac{n}{r-c}\rfloor$. In addition, let $k,d$ be non-negative integers such that $b = kt +d$ with $d < t$. We will consider two cases and exhibit an $r$-partite, $S_t$-saturated graph on $n$ vertices in each case. 
        
        Begin by supposing $a(r-c)$ is even. By rewriting (\ref{ineq3}), we see that $t-1 \leq a(r-c-1).$ 
        By Proposition \ref{hoff_const}, there exists a $(t-1)$-regular, $(r-c)$-partite graph $G$ on $a(r-c)$ vertices. Taking the disjoint union $G\cup kK_t\cup K_d$ yields an $r$-partite, $S_t$-saturated graph on $n$ vertices. 
        
        Similarly, when $a(r-c)$ is odd, there exists a $(t-1)$-regular, $(r-c)$-partite graph $G$ on $(a-1)(r-c)$ vertices. Thus $G\cup \ell K_t\cup K_m$ where $b+r-c = \ell t+m$ with $m < t$ provides an $r$-partite, $S_t$-saturated graph on $n$ vertices.
    \end{proof}
\end{theorem}
With regards to the lower bounds on $n$, we note that if $n< t+1$, then there is no $S_t$-saturated graph on $n$ vertices, and hence none that is $r$-partite. We now provide a necessary  condition for the existence of $S_t$-saturated, $r$-partite graphs that is related to the other bound in Theorem \ref{theorem_partite}.

\begin{proposition}\label{prop2.8}
    For all $r\geq 3$ and $t\geq 3$, if there exists a graph $G$ on $n$ vertices that is an $S_t$-saturated graph and $r$-partite, then $n \geq \frac{r(t-1)}{r-1}$. 
    \begin{proof}
        Suppose such a graph $G$ exists with $n < \frac{r(t-1)}{r-1}$. Let $a,b$ be non-negative integers such that $n = ar + b$ with $b < r$. Since $G$ is $r$-partite, there exists a partition $P_1,\dots,P_r$ of the vertices of $G$ with $|P_i| \leq |P_{i+1}|$ for all $1\leq i \leq r-1$ such that no two vertices in a given $P_i$ are adjacent. We have two cases to consider. 
        
        Begin by supposing that $|P_1| < |P_r|$. If $a = 0$, then $G$ has $b$ vertices with $b < r$. Since $G$ is $S_t$-saturated, we have $b\geq t+1$, and, by assumption, $b < \frac{r(t-1)}{r-1}$. It follows that $t + 1 < \frac{r(t-1)}{r-1}$. Rearranging, we have $2r < t+ 1$, a contradiction since $b < r$ and $b \geq t+1$. Thus it must be the case that $a \geq 1$.
        
        Now, since $|P_1|< |P_r|$, it follows that 
        $$|P_r| \geq \left\lceil \frac{n}{r}\right\rceil = a+1.$$ 
        For any vertex $u \in P_r$, we have $d(u) \leq n-a-1$. Since each $P_i$ is an independent set, it can contain at most one vertex of degree less than $t-1$. In particular, since $a \geq 1$, $P_r$ must contain at least one vertex of degree exactly $t-1$. Thus $n-a-1\geq t-1$. Since $n < \frac{r(t-1)}{r-1}$ by assumption, we have
        \begin{align*}
            \frac{r(t-1)}{r-1} &> t + a. \\
            \intertext{After rearranging the previous inequality and noting that $t-1\leq n-a-1 = a(r-1+b-1)$, we have }
            a(r-1) + b-1 &> a(r-1) + r-1,
        \end{align*}
        contradicting the assumption that $b < r$.
        
        Finally, we consider the case where each of the $r$ parts have equal size. It follows that $b = 0$ and $n = ar$. Every vertex in $G$ has degree at most $a(r-1)$. Since $n <\frac{r(t-1)}{r-1}$, we have that $a(r-1) < t-1$. This means that if we add an edge to $G$, the maximum degree is at most $t-1$, a contradiction to the assumption that $G$ is $S_t$-saturated.  
    \end{proof}
\end{proposition}

In addition to bridging the gap between these bounds, we would like to know which value of $c$ for given $r$ and $t$ minimizes the lower bound on $n$ in Theorem \ref{theorem_partite}. If our bound did not include ceilings, this would be a straightforward computation as demonstrated in the lemma below. However, finding a general solution using the bound in the theorem is more complicated. Although we do not have a general solution, we determine which value of $c$ provides the smallest bound on $n$ in Theorem \ref{theorem_partite} for the existence of $r$-partite, $S_t$-saturated graphs in two special cases.

\begin{lemma}\label{lemma_straightforward}
    Let $r\geq 3$ and $t\geq 3$ be fixed. Then
    $$n_2(c) = (r-c)\left(\frac{t-1}{r-c-1}\right) + r-c$$
    is minimized on the interval $[0,r-2]$ when $c = r-1-\sqrt{t-1}$.
    \begin{proof}
        Taking the derivative with respect to $c$, we obtain
        $$n_2'(c) = \frac{t-1}{(r-c-1)^2} - 1.$$
        Setting this equal to $0$, we have $c = r-1\pm\sqrt{t-1}$. Since $c\leq r-2$ and $n_2''(r-1-\sqrt{t-1}) < 0$, our function is minimized at $c = r-1-\sqrt{t-1}$. 
    \end{proof}
\end{lemma}

\begin{proposition}\label{prop_square+1}
    For all $r\geq 3$, if $t \leq (r-1)^2+1$ and $t = k^2+1$ for some integer $k$, then $r-1-\sqrt{t-1}$ is the minimizing $c$-value in Theorem \ref{theorem_partite}. The corresponding lower bound on $n$ is $n\geq (k+1)^2 = t + 2\sqrt{t-1}$.
    \begin{proof}
        Let $$n_1(c) = (r-c)\left\lceil\frac{t-1}{r-c-1}\right\rceil+r-c$$ and $$n_2(c) = (r-c)\left(\frac{t-1}{r-c-1}\right)+r-c.$$ That is, $n_1(c)$ provides the lower bounds on $n$ for the $(r-c)$-partite construction in Theorem \ref{theorem_partite}. By Lemma \ref{lemma_straightforward}, $n_2(c)$ is minimized when $c = r-1\pm \sqrt{t-1}$ on the interval $[0,r-2]$. Since $n_1(c)\geq n_2(c)$ for all $c$, it follows that $n_1(c)$ is minimized at $c = r-1-\sqrt{t-1}$ if $n_1(r-1-\sqrt{t-1}) = n_2(r-1-\sqrt{t-1})$. This is true precisely when $t =k^2+1$ for some integer $k$. We finally note that in this setting $n_1(r-1-\sqrt{t-1}) = t+2\sqrt{t-1}$.
    \end{proof}
\end{proposition}
\begin{proposition}
    For all $r\geq 3$, if $t \geq (r-1)^2 +1$ and $r-1$ divides $t-1$, then $0$ is the minimizing $c$-value in Theorem \ref{theorem_partite}. The corresponding lower bound on $n$ is $n \geq r\left(\frac{t-1}{r-1}\right)+r$.
    \begin{proof}
        Define $n_1(c)$ and $n_2(c)$ as in the previous proof. Since $t\geq (r-1)^2 + 1$, we have that $r-1-\sqrt{t-1}\leq 0$. Thus $n_2(c)$ is increasing on the interval $[0,r-2]$. It follows that $n_1(c)$ is minimized at $c = 0$, if $n_1(0) = n_2(0).$ This is true precisely when $\frac{t-1}{r-1}$ is an integer. That is, when $r-1$ divides $t-1$. We finally note that in this setting $n_1(0) = r\bigl(\frac{t-1}{r-1}\bigr)+r$.
    \end{proof}
\end{proposition}

Note that when $t > (r-1)^2+1$, Proposition \ref{prop_square+1} shows that the minimizing $c$-value for Theorem \ref{theorem_partite} is not $0$. That is, starting with fewer parts than what we're permitted leads to a smaller lower bound on $n$.


We conclude this section by providing necessary conditions for the existence of $S_t$-saturated graphs on $n$ vertices that are $K_{r+1}$-free.
\begin{lemma}\label{lemma_2.12}
    Let $G$ be an $S_t$-saturated graph on $n$ vertices with $m$ vertices of degree less than $t-1$ and no copy of $K_{r+1}$. Then
    $$n\geq \frac{r}{r-1}\left(\frac{t-1}{2}+\sqrt{\left(\frac{t-1}{2}\right)^2-\frac{m(r-1)(t-m)}{r}}\right).$$
    \begin{proof}
        Since $G$ is $K_{r+1}$-free, we have by Tur\'an's Theorem that 
        $$e(G) \leq \frac{r-1}{r}\cdot \frac{n^2}{2}.$$
        Using the structure of $S_t$-saturated graphs described in Lemma \ref{lemma_2.1} at the beginning of this section, we also have that 
        $$e(G) \geq \frac{1}{2}(n-m)(t-1) + {{m}\choose{2}}$$ for some value of $m$. These two bounds together provide the desired inequality.
    \end{proof}
\end{lemma}

Considering the case where $t\geq 2r$, we provide a general bound independent of the existence of $S_t$-saturated, $K_r$-free graphs with specific $m$. We also demonstrate when these bounds are at least $t+1$, the trivial necessary condition for the existence of an $S_t$-saturated graph.

\begin{proposition}
    Let $G$ be an $S_t$-saturated graph on $n$ vertices with no copy of $K_{r+1}$ where $t\geq 2r$ and $r\geq 2$. Then 
    $$n \geq \frac{r}{r-1}\left(\frac{t-1}{2}+\sqrt{\left(\frac{t-1}{2}\right)^2-(r-1)(t-r)}\right).$$
    Furthermore, this bound is at least $t + 1$ whenever $t \geq r(r+1)-1$. 
    \begin{proof}
        We begin by noting that the bound in Lemma \ref{lemma_2.12} is minimized at $m = r$. This is because $m\leq r$ by Lemma \ref{lemma_2.1} and because the bound in Lemma \ref{lemma_2.12} is decreasing in $m$ for $m < t/2$. Plugging this value in for $m$ yields the desired bound. Letting $n(t)$ denote this bound, we observe that $n(r(r+1)-1) = r(r+1)=t+1$. In addition, $n'(t)\geq \frac{r}{r-1}$ for all $t$. By the racetrack principle, $n(t) \geq t +1$ whenever $t\geq r(r+1)-1$ as desired.
    \end{proof}
\end{proposition}

Finally, we observe that the bound in Lemma \ref{lemma_2.12} is asymptotically equivalent to the bound in Proposition \ref{prop2.8}.
\begin{proposition}
    Let $r$ be fixed. If there exists an $S_t$-saturated graph on $n$ vertices with no copy of $K_{r+1}$ where $t\geq 2r$ and $r\geq 3$, then $$n \geq \frac{r(t-1)}{r-1} - O_t(1).$$
    \begin{proof}
        Note that the number of vertices, that is $m$, of degree less than $t-1$, in a $K_{r+1}$-free graph is at most $r$ since these vertices form in a clique in $S_t$-saturated graphs. Thus the following holds as $t \geq 2r$.
    \begin{align*}
        \frac{r}{r-1}\left(\frac{t-1}{2}+\sqrt{\left(\frac{t-1}{2}\right)^2-\frac{m(r-1)(t-m)}{r}}\right) 
        &\ge  \frac{r}{r-1}\left( \frac{t-1}{2} + \frac{t-1}2 - \frac{m(r-1)t}{r(t-1)}  \right) \\
        &= \frac{r}{r-1}\left(t-1-O_t(1) \right) \\
        &= \frac{r(t-1)}{r-1}-O_t(1). \qedhere
    \end{align*}
    \end{proof}
\end{proposition}

\section{Stars in star-saturated graphs}\label{sec_star_in_star}
\subsection{General results for stars in star-saturated graphs}\label{sec_starstar1}
We now turn to counting copies of stars $S_r$ in $S_t$-saturated graphs. As stated in the introduction, we write $s_r(G)$ to denote the number of copies of $S_r$ in a given graph $G$ where $S_r$ is the complete bipartite graph $K_{1,r}$. Note that if $t \leq r$, then $\sat_{S_r}(n,S_t) = 0$ trivially as an $S_t$-saturated graph must be $S_t$-free. That is, any $S_t$-saturated graph has no vertex of degree at least $t$, and hence none of degree at least $r$. Our focus is therefore on the situation where $t > r$, and we will consider the cases where $t$ is odd and even separately. We also ignore the case where $t= 1$ as the only $S_1$-saturated graph on $n$ vertices is an independent set. 

We stated in Section \ref{Constructions} that $\KR_{t,n}(m)$ would be a candidate for achieving $\sat_{S_r}(n,S_t)$. Our next theorem shows that this graph does exactly that.

\begin{theorem}\label{theorem_3.2}
    For all $n\geq 2t-1$ with $t \geq 2$ and $r< t$,  
    \[ \sat_{S_r}(n,S_t) =  
        \min_{0\leq m \leq t-1}
        s_r(\KR_{t,n}(m)).
    \]
    Note also that 
    \[
        s_r(\KR_{t,n}(m)) = 
            \begin{cases}
                m{{m-1}\choose{r}} + (n-m){{t-1}\choose{r}} & \text{if $(t-1)(n-m)$ is even} \\
                m{{m-1}\choose{r}} + (n-m){{t-1}\choose{r}} + {{m-1}\choose{r-1}} & \text{if $(t-1)(n-m)$ is odd.}
            \end{cases}
    \]
    \begin{proof}
        We begin by considering the case where $t$ is odd. 
        By Lemma \ref{lemma_2.1}, if $G$ is $S_t$-saturated, then $G$ contains a clique $A$ containing all of the vertices with degree smaller than $t-1$. Let $A$ have size $m$, and let $B=V(G)\setminus A$. We have two cases to consider. If $G$ has no edges between $A$ and $B$, then $G$ contains exactly $m{{m-1}\choose{r}} + (n-m){{t-1}\choose{r}}$ copies of $S_r$. The first term counts stars centered in $A$, and the second term counts stars centered in $B$. Since $t-1$ is even, there exists a $(t-1)$-regular graph $R_{t-1,n-m}$ on $n-m$ vertices for all $n-m\geq t$. This inequality holds since we assume $n\geq 2t-1$ and $m \leq t-1$. Thus an $S_t$-saturated graph with precisely the above count is given by $\KR_{t,n}(m)$. 
        
        Now, if there exist vertices $u\in A$ and $v\in B$ such that $u$ is adjacent to $v$, then our graph contains all of the previously counted stars, along with at least ${{m-1}\choose{r-1}}$ stars centered at $u$ containing the edge $uv$. This means than an $S_t$-saturated graph $G$ with $|A| = m$ and any such edge must have at least as many copies of $S_r$ as $\KR_{t,n}(m)$. Thus an $S_t$-saturated graph with minimum number of copies of $S_r$ is given by some $\KR_{t,n}(m)$ for some $m\leq t-1$. 
        
        We now consider the case where $t$ is even. For a given $m \leq t-1$, if $n-m$ is even, we can find a $(t-1)$-regular graph on $n-m$ vertices, 
        and the argument is the same as before. That is, among $S_t$-saturated graphs with $m$ vertices of degree less than $t-1$, $\KR_{t,n}(m)$ is a minimal example with respect to copies of $S_r$. When $n-m$ is odd, we can construct a graph $R_{t-1,n-m}$ on $n-m$ vertices that is $(t-1)$-regular with the exception of one vertex $v$ of degree $t-2$. Thus we can construct an $S_t$-saturated graph $\KR_{t,n}(m)$ by taking the disjoint union of $K_m$ with $R_{t-1,n-m}$ and adding an edge from $v$ to an arbitrary vertex in the clique $K_m$. Since every $S_t$-saturated graph with $m$ vertices of degree less than $t-1$ has at least $m{{m-1}\choose{r}} + (n-m){{t-1}\choose{r}}$ many copies of $S_r$ and there is no $(t-1)$-regular graph on $n-m$ vertices, this adds the fewest possible copies of $S_r$. That is, we must have at least one edge between $A$ and $B$, introducing ${{m-1}\choose{r-1}}$ copies of $S_r$. Therefore the generalized saturation number is obtained by minimizing $s_r(\KR_{t,n}(m))$ over all values of $m$ between the two scenarios.
    \end{proof}
\end{theorem}

\begin{figure}[H]
    \centering
    \begin{tikzpicture}
        \draw[fill=black] (0,0) circle (3pt);
        \draw[fill=black] (0,1) circle (3pt);
        \draw[fill=black] (1,2) circle (3pt);
        \draw[fill=black] (1,-1) circle (3pt);
        \draw[fill=black] (1.5,1) circle (3pt);
        \draw[fill=black] (1.5,0) circle (3pt);
        \draw[thick] (1.5,1) -- (1,2) -- (1,-1) -- (1.5,1) -- (1.5,0) -- (1,2) -- (1,-1) -- (1.5,0);
        \draw[thick] (1,2) -- (0,1) -- (0,0) -- (1,-1);
        \draw[thick] (0,1) -- (0,0);
        \draw[thick] (0,1) -- (1.5,1);
        \draw[thick] (0,0) -- (1.5,0);
    \end{tikzpicture}
    \caption{An example of an $S_5$-saturated graph that does not satisfy the criteria in Theorem \ref{theorem_3.2}}
    \label{Fig_3.2_counter}
\end{figure}
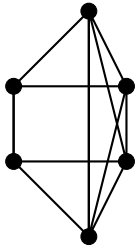

We note that the above theorem does not hold when $n < 2t-1$. This is because it is possible for an $S_t$-saturated graph to have fewer than $t$ vertices of degree $t-1$. This means we can't consider the disjoint union of a small clique and a $(t-1)$-regular graph. In this setting, it turns out than an optimal graph does not need to have that structure. For example, a quick check shows that the graph in Figure \ref{Fig_3.2_counter} minimizes the number of copies of $S_3$ among $S_5$-saturated graphs on $6$ vertices and has only $4$ vertices of degree $4$. 

We now introduce additional notation to aid in our discussion of this topic. Given $n,r,$ and $t$ with $n \geq \max\{2t-1,t+1\}$, we define 
\[m_0(n,r,t) := \underset{m}{\argmin}\  s_r(\KR_{t,n}(m)).
\]
That is, $m_0(n,r,t)$ is the value of $m$ for which $\KR_{t,n}(m)$ attains the generalized saturation number $\sat_{S_r}(n,S_t)$. If the generalized saturation number is achieved for multiple values of $m$, we take the smallest one for definiteness. In light of Theorem \ref{theorem_3.2}, our goal is to identify the value of $m_0(n,r,t)$ for given values of $n,r,$ and $t$. Kaszonyi and Tuza \cite{KT} showed that the number of edges in an $S_t$-saturated graph is minimized when $m =\lfloor \frac{t}{2}\rfloor$ or $\lfloor \frac{t+1}{2} \rfloor$, answering our question for $r = 1$. This solution does not hold for all $r\geq 1$ though. Rather $m_0(n,r,t)$ depends on both $r$ and $t$. The value of $n$ does not matter for $n \geq 2t-1$ since we are simply getting more vertices of degree $t-1$ as we increase $n$. The particular value of $n$ does matter when $t+1 \leq n < 2t-1$ though. This is because there may not exist a $(t-1)$-regular graph on $n-m$ as demonstrated by the example in Figure \ref{Fig_3.2_counter}.

Although we are unable to provide a closed form for $m_0(n,r,t)$ for arbitrary pairs of $r$ and $t$, there is more we can say about the optimal choice, or in some cases, choices. We begin with the observation that 
\begin{equation}\label{eq1}
    \mathcal{D}(m) := s_r(\KR_{t,n}(m+1)) - s_r(\KR_{t,n}(m)) = (r+1){{m}\choose{r}} - {{t-1}\choose{r}}.
\end{equation}
That is, $\mathcal{D}(m)$ denotes the change in the number of copies of $S_r$ as we increase $m$, the number of vertices of degree less than $t-1$, by $1$. The following is an immediate consequence of this observation. 

\begin{lemma}
    If $t$ is odd, then for all $n\geq 2t-1$ with $t > r$, we have $m_0(n,r,t) \geq r$. In particular, $m_0(n,t-1,t) = t-1$.  
    \begin{proof}
        If $m < r$, then ${{m}\choose{r}} = 0$ and so $\mathcal{D}(m) < 0$. That is, $s_r(\KR_{t,n}(m+1)) < s_r(\KR_{t,n}(m))$ whenever $m < r$. Therefore the minimum must be attained when $m\geq r$. To prove the second statement, we observe that when $r = t-1$, $\mathcal{D}(t-1) = t-1 > 0$. Since $\mathcal{D}$ is an increasing function, it follows that $m_0(n,t-1,t) = t-1$.
    \end{proof}
\end{lemma}

For fixed $t$ and $r$, we extend the definition of $\mathcal{D}$ in (\ref{eq1}) to all real numbers: 
    $$\mathcal{D}(x) := (r+1){{x}\choose{r}} - {{t-1}\choose{r}}.$$
Now, ${{x}\choose{r}}$ is convex and increasing in $x$ for all $x\geq r- 1$. It follows that those properties hold for $\mathcal{D}(x)$ as well. Thus $\mathcal{D}(x)$ has a unique root on the interval $(r-1,\infty)$ as $\mathcal{D}(r-1) = -{{t-1}\choose{r}}$. Our next theorem takes advantage of this structure.

\begin{theorem}\label{theorem_3.5}
    For fixed $n,r,$ and $t$ with $n\geq 2t-1$ and $t>r$, let $\overline{x}$ denote the unique root of $\mathcal{D}(x)$ in the interval $(r-1,\infty)$. Then 
    $$m_0(n,r,t) = \lceil \overline{x}\rceil.$$
    Furthermore, when $\overline{x} \not\in \Z$, this is the unique minimizing value of $m$ in Theorem $\ref{theorem_3.2}$. When $\overline{x} \in \Z$, both $\overline{x}$ and $\overline{x} + 1$ simultaneously minimize the number of copies of $S_r$ among $S_t$-saturated graphs on $ n$ vertices.
    \begin{proof}
        As stated previously, $\mathcal{D}(x)$ has a unique root $\overline{x}$ and is increasing on the interval $[r-1,\infty)$. Let $m\geq r$ be an integer. If $m < \overline{x}$, then $\mathcal{D}(m) < 0$ and $s_r(\KR_{t,n}(m+1)) < s_r(\KR_{t,n}(m))$. If $m > \overline{x}$, then $\mathcal{D}(m) > 0$ and $s_r(\KR_{t,n}(m+1)) > s_r(\KR_{t,n}(m))$. If $\overline{x} \not\in \Z$, then it follows that $s_r(\KR_{t,n}(m))$ is minimized when $m$ is the first integer larger than $\overline{x}$, namely $\lceil \overline{x} \rceil$, and this choice of $m$ is unique. 
        
        On the other hand, if $\overline{x}\in \Z$, then $\mathcal{D}(\overline{x}) = 0$ and $s_r(\KR_{t,n}(\overline{x}+1)) = s_r(\KR_{t,n}(\overline{x}))$. Therefore $s_r(\KR_{t,n}(m))$ is minimized by $\overline{x}$ and $\overline{x}+1$ simultaneously. These are the only optimal choices for $m$ as $\mathcal{D}(\overline{x}-1) < 0$ and $\mathcal{D}(\overline{x}+1) > 0$.
    \end{proof}
\end{theorem}

We conclude the discussion on general results by providing two lower bounds on the value of $m_0(n,r,t)$ and giving a more precise answer for $t = 2$. The second lower bound will be of additional interest in the following section on asymptotic results. Before stating these results, we prove a simple lemma regarding binomial coefficients. 

\begin{lemma}\label{bin_lemma1}
    If $a\geq c \geq 2$ and $b> 1$ where $a,b\in \R$ and $c \in \Z$, then $$b^c{{\lfloor a/b\rfloor}\choose{c}} < {{a}\choose{c}}.$$
    \begin{proof}
        Note that if $\lfloor a/b\rfloor < c$, then the inequality holds trivially as the left hand side of our inequality is equal to $0$, and the right hand side is positive. 
        Suppose then that $\lfloor a/b\rfloor \geq c$. Then we have the following
        \begin{align*}
            b^c(\lfloor a/b\rfloor)_c &= b^c(\lfloor a/b\rfloor)(\lfloor a/b\rfloor-1)\cdots(\lfloor a/b\rfloor-(c-1)) \leq  a(a-b)(a-2b)\cdots(a-b(c-1)) \\ 
            \intertext{and}
            (a)_c &= a(a-1)(a-2)\cdots(a-(c-1)).
        \end{align*}
        Since $0 < a-bk < a - k$ for all $1 \leq k\leq c-1$ and $b > 1$, the desired inequality holds.
    \end{proof}
\end{lemma}

We now proceed to state and prove our lower bounds on the optimal choice for $m$.

\begin{corollary}\label{cor_3.8}
    If $t\geq 3$ is odd and $r \geq 2$ with $t> r$, then for all $n\geq 2t-1$, we have $$m_0(n,r,t) \geq \frac{t+1}{2}.$$
    \begin{proof}
        Since $r \geq 2$, we know that $2^r > r + 1$. Let $m\leq  \frac{t-1}{2}$ be an integer. Applying Lemma \ref{bin_lemma1} with $a = t-1,b = 2$, and $c = r$, we obtain the following
        $$(r+1){{m}\choose{r}} < 2^r{{m}\choose{r}} \leq 2^r{{\frac{t-1}{2}}\choose{r}} < {{t-1}\choose{r}}.$$
        Thus $$\mathcal{D}(m)= (r+1){{m}\choose{r}} - {{t-1}\choose{r}} < 0.$$
        This means that $\sat_{S_r}(n,S_t)$ is not attained by $\KR_{t,n}(m)$, and $m_0(n,r,t) > \frac{t-1}{2}$. 
    \end{proof}
\end{corollary}

\begin{corollary}
    If $t\geq 3$ is odd and $r\geq 2$ with $t > r$, then for all $n\geq 2t-1$, we have $$m_0(n,r,t) > \frac{t-1}{(r+1)^{1/r}}.$$ 
    \begin{proof}
        Note that $$\mathcal{D}\left(\left\lfloor \frac{t-1}{(r+1)^{1/r}}\right\rfloor \right) = (r+1){{\left\lfloor \frac{t-1}{(r+1)^{1/r}}\right\rfloor}\choose{r}} - {{t-1}\choose{r}}.$$
        Applying Lemma \ref{bin_lemma1} with $a = t-1$, $b = (r+1)^{1/r}$, and $c = r$, we see that this quantity is strictly less than $0$. Since $\frac{t-1}{(r+1)^{1/r}}$ is not an integer for all $r\geq 2$, it must be the case that $m_0(n,r,t) > \frac{t-1}{(r+1)^{1/r}}$.
    \end{proof}
\end{corollary}

When $r = 2$, finding $m_0(n,r,t)$ amounts to solving a quadratic equation and applying Theorem \ref{theorem_3.5}. The following is thus immediate. 

\begin{proposition}\label{prop_3.10}
    For all $t \geq 3$ and $n\geq 2t-1$, the value of $\overline{x}$ as in Theorem \ref{theorem_3.5} for $r = 2$ is given by 
    $$\overline{x} = \frac{1}{2} +\frac{1}{6}\sqrt{12t^2-36t+33}.$$
\end{proposition}

We can say a little more when $r =2$.
\begin{proposition}
    There are two optimal choices for $m$ that simultaneously minimize the number of copies of $S_2$ among $S_t$-saturated graphs if and only if $t$ is given by the following where $i \geq 0$ is some non-negative integer
    $$t(i) = \frac{1}{4}\left((1+\sqrt{3})(2+\sqrt{3})^{i}-(\sqrt{3}-1)(2-\sqrt{3})^{i}-2\right)+2.$$
    \begin{proof}
        By Theorem \ref{theorem_3.5} and Proposition \ref{prop_3.10}, there are two optimal choices for $m$ precisely when $m = \frac{1}{2}+\frac{1}{6}\sqrt{12t^2-36t+33}$ is an integer. This is the case when we can write $\sqrt{12t^2-36t+33}$ in the form $6k + 3$ where $k$ is an integer. Equivalently, we need $(t-1)(t-2) = 3(k+1)k$. 
        Now, the second member of the Diophantine pair $(x,y)$ that satisfies $3(x^2 + x) = y^2 + y$ is given by $y = a(i)$ where $a(i)$ satisfies the recurrence (see OEIS sequence A001571 \cite{OEIS}) $$a(i) = 4a(i-1) - a(i-2) + 1 \text{ with } a(0) = 0 \text{ and } a(1)=2.$$ Solving the linear recurrence, we find that $$a(i) = \frac{1}{4}\left((1+\sqrt{3})(2+\sqrt{3})^{i}-(\sqrt{3}-1)(2-\sqrt{3})^{i}-2\right).$$
        Thus the values of $t$ for which there are two optimal choices of $m$ are given by $t(i) = a(i) + 2$. 
    \end{proof}
\end{proposition}

\subsection{Asymptotic results for stars in star-saturated graphs}\label{sec_star_asymp}
In addition to general questions regarding the number of stars in star-saturated graphs, we can address asymptotic questions. Here we focus on the case where $t$ is odd and $n\geq 2t-1$ for convenience. Utilizing our general results from the previous section, we immediately proceed to our asymptotic results.

\begin{theorem}
    Let $r = o(\sqrt{t})$ with $t$ odd and $r\geq 2$. Then for all $n\geq 2t-1$, we have $$m_0(n,r,t) = (1+o_t(1))\frac{t-1}{(r+1)^{1/r}}.$$ 
    \begin{proof}
        As in Theorem \ref{theorem_3.5}, we need to solve for the value of $x$ such that
        \begin{equation}\label{eq3}
            (r+1){{x}\choose{r}} = {{t-1}\choose{r}}.
        \end{equation}
        Let $\overline{x}$ be the unique solution. 
        By Corollary \ref{cor_3.8}, we know that $\overline{x}\geq \lfloor t/2\rfloor$ and so $\overline{x} = \Theta(t)$. It is well known that when $k = o(\sqrt{n})$, $${{n}\choose{k}} = (1+o_n(1))\frac{n^k}{k!}.$$  Since $r =o(\sqrt{t})$, $\overline{x} = \Theta(t)$, and $n\geq 2t-1$, we can apply this to our equality and obtain the following.
        \begin{align*}
            (r+1)(1+o_{\overline{x}}(1))\overline{x}^r &= (1+o_t(1))(t-1)^r. \\
            \intertext{Thus, since $\overline{x} = \Theta(t)$,} 
            \overline{x} &= (1+o_t(1))\frac{t-1}{(r+1)^{1/r}}.
        \end{align*} Applying Theorem \ref{theorem_3.5} gives us the desired result.
    \end{proof}
\end{theorem}

With slightly less precision than our previous theorem, we consider the more general case where $r = o(t)$. Before stating our theorem, we note the following useful result. A short proof can be found in \cite{D}. See also \cite{M}.

\begin{lemma}\label{lemma_as}
    If $k = o(n)$, then $\log{{n}\choose{k}} = (1+o_n(1))k\log\frac{n}{k}$.
\end{lemma}

\begin{theorem}
    Let $r = o(t)$ with $t$ odd and $r\geq 2$. Then for all $n\geq 2t-1$, we have that $$m_0(n,r,t) = \left(\frac{t-1}{(r+1)^{1/r}}\right)^{1+o_t(1)}.$$
    \begin{proof}
        Let $\overline{x}$ be the unique solution to (\ref{eq3}). By Corollary \ref{cor_3.8}, $\overline{x} \geq \lfloor t/2 \rfloor$ and so $\overline{x} = \Theta(t)$. We now take the logarithm of both sides in our equality, and we observe the following. 
        \begin{align*}
            \log\left((r+1){{\overline{x}}\choose{r}}\right) &= \log{{t-1}\choose{r}}, \text{ so} \\
            \log(r+1) + \log{{\overline{x}}\choose{r}} &= \log{{t-1}\choose{r}} \\
            \intertext{By Lemma \ref{lemma_as},}
            \log(r+1) + (1+o_{\overline{x}}(1))r \log\left(\frac{\overline{x}}{r}\right) &= (1+o_t(1))r\log\left(\frac{t-1}{r}\right) 
        \end{align*}
        Solving for $\overline{x}$ in this equation and using the fact that $\overline{x} =\Theta(t)$, we get
        $$\overline{x} = \left(\frac{t-1}{(r+1)^{1/r}}\right)^{1+o_t(1)}.$$ 
        By Theorem \ref{theorem_3.5}, $m_0(n,r,t) = \lceil \overline{x}\rceil$, and the result follows.
    \end{proof}
\end{theorem}

\section{Stars in clique-saturated graphs}\label{sec_star_clique}
Having considered the minimum number of cliques of a given size in star-saturated graphs and counting stars in star-saturated graphs, we move to counting stars in clique-saturated graphs. Clique-saturated graphs have been of interest since the proof of Tur\'an's Theorem, which established the maximum number of edges among clique-saturated graphs. Theorem \ref{theorem_ehm} gives the minimum number of edges among clique-saturated graphs. This gives us the exact value of $\sat_{S_1}(n,K_t)$ since $S_1$ is simply an edge. The split graph $K_{t-2}+\overline{K}_{n-t+2}$ satisfies the edge count for Theorem \ref{theorem_ehm}. For $n$ sufficiently large, Chakraborti and Loh \cite{CL} showed that this graph is the unique $K_t$-saturated graph minimizing the number of copies of $K_r$ for all $r < t$, as well as the number of copies of $C_r$.

However, it turns out that the split graph is far from optimal when we want to minimize stars $S_r$ for $r\geq 3$. Using a construction of Alon, Erd\H{o}s, Holzman, and Krivelevich \cite{AEHM} involving a truncated projective plane for their upper bound, Ergemlidze, Methuku, Tait, and Timmons \cite{EMTT} proved the following.

\begin{theorem}[Ergemlidze et al., 2021]
    For integers $n\geq t\geq 3$ and $r\geq 3$,
    $$\sat_{S_r}(n,K_t) = \Theta(n^{r/2+1})$$
\end{theorem}

An important detail about the construction from Alon, Erd\H{o}s, Holzman, and Krivelevich is that the resulting graph has maximum degree $\Delta = \Theta(\sqrt{n})$. Since vertices of large degree contain many stars, it makes sense that the optimal choice is for the graph to have as small of degree as possible. Furthermore, Alon et al. showed that the maximum degree of a $K_t$-saturated graph on $n$ vertices is $\Omega(\sqrt{n})$ \cite{AEHM}.

In addition to the split graph not minimizing stars among all $K_t$-saturated graphs, it is still not optimal for minimizing stars when considering families of $K_t$-saturated graphs with linear maximum degree. The following graph will be the starting point for our construction.

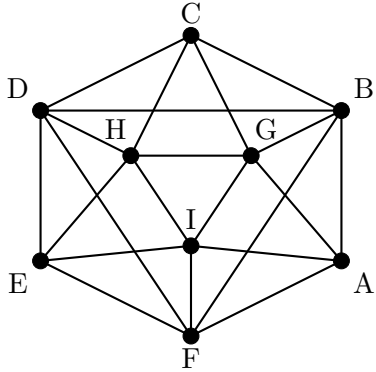
\begin{figure}[H]
    \centering
    \begin{tikzpicture}
        \draw[fill=black] (2,-1) circle (3pt);
        \draw[fill=black] (2,1) circle (3pt);
        \draw[fill=black] (0,2) circle (3pt);
        \draw[fill=black] (-2,1) circle (3pt);
        \draw[fill=black] (-2,-1) circle (3pt);
        \draw[fill=black] (0,-2) circle (3pt);
        \draw[fill=black] (0.8,0.4) circle (3pt);
        \draw[fill=black] (-0.8,0.4) circle (3pt);
        \draw[fill=black] (0,-0.8) circle (3pt);
        \node at (2.3,-1.3) {A};
        \node at (2.3,1.3) {B};
        \node at (0,2.3) {C};
        \node at (-2.3,1.3) {D};
        \node at (-2.3,-1.3) {E};
        \node at (0,-2.3) {F};
        \node at (1,0.75) {G};
        \node at (-1,0.75) {H};
        \node at (0,-0.45) {I};
        \draw[thick] (2,-1) -- (2,1) -- (0,2) -- (-2,1) -- (-2,-1) -- (0,-2) -- (2,-1);
        \draw[thick] (2,1) -- (-2,1) -- (0,-2) -- (2,1);
        \draw[thick] (0.8,0.4) -- (-0.8,0.4) -- (0,-0.8) -- (0.8,0.4);
        \draw[thick] (-0.8,0.4) -- (-2,-1) -- (0,-0.8) -- (2,-1) -- (0.8,0.4) -- (0,2) -- (-0.8,0.4);
        \draw[thick] (-2,1) -- (-0.8,0.4);
        \draw[thick] (0,-2) -- (0,-0.8);
        \draw[thick] (2,1) -- (0.8,0.4);
    \end{tikzpicture}
    \caption{$K_4$-saturated graph $G_{4,9}$ on $9$ vertices}
    \label{fig_sat}
\end{figure}

\begin{proposition}
    Let $t\geq 4$ and $r\geq 3$. 
    There exists a sequence $(G_{t,n})$ of $K_t$-saturated graphs on $n$ vertices with $\Delta(G_{t,n}) = \Theta(n)$ and a constant $n_{r,t}$ such that $s_r(G_{t,n}) < s_r(K_{t-2}+\overline{K}_{n-t+2})$ for all $n\geq n_{r,t}$. 
    \begin{proof}
        Consider the graph $G_{4,9}$ in Figure \ref{fig_sat}. We obtain $G_{4,n}$ for $n > 9$ by blowing up vertices $A,C,E$ into independent sets of size as equal as possible. For $t > 4$, we define $G_{t,n}$ to be $G_{4,n-t+4}+K_{t-4}$. Since $G_{4,n-t+4}$ is $K_4$-saturated, joining $t-4$ universal vertices results in a $K_t$-saturated graph. Given $r$ and $t$, we can find a constant $n_{r,t}$ such that $s_r(G_{t,n})< s_r(K_{t-2}+\overline{K}_{n-t+2})$ for all $n\geq n_{r,t}$.
    \end{proof}
\end{proposition}
Although we omit the computations from the previous proof, we note that the key to the construction is that $\Delta(G_{t,n})$ is roughly $\frac23n$ for large $n$, and $\Delta(K_{t-2}+\overline{K}_{n-t+2}) = n-1$.

\section{Paths and cliques}\label{sec_paths}
In addition to stars and cliques, we are interested in the broader topic of saturation involving trees and cliques. In the next section we will say a few things about arbitrary trees, but here we will turn to another basic class of trees, namely paths. Let $P_{r+1}$ denote a path on $r+1$ vertices. We begin this discussion by stating the following result of Kritschgau, Methuku, Tait, and Timmons \cite{KMTT} concerning cycles.
\begin{proposition}[Kritschgau et al., 2020]
    For $t\geq 5$ and $r\leq 2t-4$, $$\sat_{C_r}(n,K_t) = \Theta(n^{\lfloor r/2\rfloor}).$$
\end{proposition}

The case where $t = 4$ was accounted for in their result on cliques in clique-saturated graphs. By refining their proof slightly and with a similar argument used by Chakraborti and Loh \cite{CL}, we prove the following  regarding the order of $n$ when counting paths in clique-saturated graphs.
\begin{theorem}
    For $t\geq 4$ and $r\leq 2t-3$, $$\sat_{P_{r+1}}(n,K_t) = \Theta(n^{\left\lceil \frac{r+1}{2}\right\rceil}).$$ If $r \geq 2t-2$, the split graph is $P_{r+1}$-free.
    \begin{proof}
        By considering the split graph $K_{t-2}+\overline{K}_{n-t+2}$, we can construct a path on $r+1$ vertices by using at most $k = \lceil \frac{r+1}{2}\rceil$ vertices from the independent set. The remaining vertices must come from the clique of order $t-2$. Thus the number of copies of $P_{r+1}$ in $G$ is $$\frac{1}{2}{{n-t+2}\choose{k}}(t-2)_kk! +o(n^{k}).$$ The $o(n^{k})$ term accounts for any paths using fewer than $k$ elements in the independent set. This gives us the appropriate upper bound. 
        
        For the lower bound, let $G$ be $K_t$-saturated, and let $I$ be an independent set of order $k =\lceil \frac{r+1}{2}\rceil$ in $G$. There are $k!$ ways to order the elements of $I$. Enumerate the elements $v_1,v_2,\dots,v_k$. For each ordering we will give a lower bound on the number of copies of $P_{r+1}$ containing it. For all $1\leq i \leq k-1$, let $V_i$ be a set of vertices such that $V_i\subseteq N(v_i)\cap N(v_{i+1})$ and the subgraph induced by $V_i$ is a copy of $K_{t-2}$. 
        Such copies exist since $v_i,v_{i+1} \in I$ which is an independent set, and $G$ is $K_t$-saturated. Since each $V_i$ has $t-2$ elements, we can pick distinct $u_i \in V_i$ such that $v_1u_1v_2\cdots v_{k-1}u_{k-1}v_k$ is a path in $G$. This gives us at least $\frac{1}{2}k!(t-2)_k$ copies of $P_{r+1}$ involving every element of $I$. The factor of $\frac{1}{2}$ accounts for the double-counting of paths being read from left-to-right and right-to-left. Chakraborti and Loh showed in \cite{CL} that any $K_t$-saturated graph contains $\Theta(n^k)$ independent sets of order $k$ for any given $k$. This gives us the corresponding lower bound, and we have that $\sat_{P_{r+1}}(n,K_t) = \Theta(n^k)$ as desired. 
        
        For the second statement in the theorem, we note that when $r\geq 2t-2$, a copy of $P_{r+1}$ in $K_{t-2}+\overline{K}_{n-t+2}$ must use at least $\lfloor\frac{r+1}{2}\rfloor$ vertices from the clique $K_{t-2}$. Thus $$2t-2 = 2(t-2)+2 \geq 2\left\lfloor\frac{r+1}{2}\right\rfloor +2 > r,$$ a contradiction. Therefore $K_{t-2}+\overline{K}_{n-t+2}$ is $P_{r+1}$-free.
    \end{proof}
\end{theorem}
Chakraborti and Loh \cite{CL} showed that for $n$ sufficiently large in terms of $r$ and $t$, the split graph minimizes the number of copies of $C_r$ in $K_t$-saturated graphs. This leads us to the following question.

\begin{question}
    For $n$ sufficiently large, is the number of copies of $P_{r+1}$ in $K_t$-saturated graphs is minimized by $K_{t-2}+\overline{K}_{n-t+2}$?
\end{question}

Shifting perspectives, we now briefly consider path-saturated graphs. Kaszonyi and Tuza provided a construction in \cite{KT} for trees which are $P_{t+1}$-saturated for all $t\geq 3$. Furthermore, they characterized all trees which are $P_{t+1}$-saturated for $t\geq 4$. The existence of such graphs, which are of course $K_3$-free, solves the asymptotic question for minimizing cliques in path-saturated graphs immediately. The lower bound on $n$ comes from \cite{KT}.
\begin{proposition}
    If $t\geq 3$ is odd, then $\sat_{K_r}(n,P_{t+1}) = 0$ for all $n\geq 3\cdot 2^{\frac{t+1}{2}-1}-2$ and all $r\geq 3$. If $t\geq 4$ is even, then $\sat_{K_r}(n,P_{t+1}) =0$ for all $n\geq 2^{\frac{t}{2}+1}-2$ and all $r\geq 3$.
\end{proposition}

\section{More general cases}\label{sec_general}
Counting paths and stars in clique-saturated graphs, as well as the reverse, is part of the larger question of counting trees in clique-saturated graphs and counting cliques in tree-saturated graphs. Some results for the traditional saturation problem on trees were proved by Kaszonyi and Tuza \cite{KT}. Many additional results on this subject are due to Faudree, Faudree, Gold, and Jacobson \cite{FFGJ}. Here we prove some initial results regarding the generalized saturation problem for trees and cliques. We let $n_T(G)$ denote the number of copies of a tree $T$ in $G$. 
\begin{proposition}
    Let $G$ be a $K_t$-saturated graph with minimum degree $\delta$, and let $r \leq 2t-4$. Then for any tree $T$ on $r$ vertices with minimum degree $\delta$, the following hold
    \begin{align}
        n_T(G) &> 0 \text{ if } \delta = t-2 \\
        n_T(G) &> 0 \text{ if } \delta = t-1 \\
        n_T(G) &> 0 \text{ if } \delta \geq 2t-5.
    \end{align}
    \begin{proof}
        Our aim is to show that for any given $T$ and $K_t$-saturated graph $G$, these minimum degree restrictions guarantee at least one copy of $T$ in $G$. Kritschgau et al. \cite{KMTT} demonstrated the exact structure of $K_t$-saturated graphs with minimum degree $t-2$ and $t-1$. In particular, when $\delta = t- 2$, $G$ is the split graph $K_{t-2}+\overline{K}_{n-t+2}$. When $\delta = t-1$, $G$ is isomorphic to $(K_{t-1}-e) +\overline{K}_{n-t+1}$ or $W_t(m_1,1,m_3,m_4,1)$ for some $m_1+m_3 + m_4 =n-t+1$. Here $e$ is any edge in $K_{t-1}$, and $W_t(m_1,m_2,m_3,m_4,m_5)$ is the graph obtained by taking a wheel with five vertices on the outer cycle and replacing the central vertex with a clique of size $t-3$ and each vertex $v_i$ of the outer cycle with an independent set of size $m_i$. Two vertices are adjacent in $W_t(m_1,m_2,m_3,m_4,m_5)$ if and only if they replaced adjacent vertices in the original wheel. We now note that for the above mentioned graphs, we can utilize the large clique to help find any tree $T$ on $t$ vertices. For the final case, we note the well known result that if a graph has minimum degree $\delta$, then it must contain any tree on $\delta+1$ vertices.
    \end{proof}
\end{proposition}

Mimicking a technique of Kaszonyi and Tuza, we provide a lower bound on $r$ for which  $\sat_{K_r}(n,F) =0$ for a $t$-vertex graph $F$ in terms of its independence number.  To this end, let $u(F) = t-\alpha(F)-1$ where $\alpha$ is the independence number of $F$. Let $d$ be the minimum number of edges in a subgraph of $G$ induced by an independent set $S$ of size $\alpha$ and one other vertex $v$. Note that the graph induced by $v$ and $S$ is the star $S_d$ and some number of isolated vertices. It is trivially true that an $F$-saturated graph can not contain any copy of $K_t$ since $F$ has $t$ vertices. The following proposition shows that we can find $F$-saturated graphs whose largest cliques are smaller. Before proving our result, we state a key lemma of Kaszonyi and Tuza \cite{KT}. 

\begin{lemma}[Kaszonyi and Tuza, 1986] Let $\mathcal{F}' = \{F_i\setminus \{x\}:x \in V(F_i), F_i \in \mathcal{F}\}$ and suppose that some vertex $x \in V(G)$ has degree $d(x) = n-1$. Then $G$ is $\mathcal{F}$-saturated if and only if $G\setminus \{x\}$ is $\mathcal{F}'$-saturated.
\end{lemma}

Here we say that a graph $G$ is $\mathcal{F}$-saturated for a family $\mathcal{F}$ of forbidden subgraphs $F_1,\dots,F_k$ if $G$ contains no $F_i$ but the addition of any missing edge creates at least one copy of some $F_i$.

\begin{proposition}
    Let $F$ be a graph on $t$ vertices. Then for $n$ sufficiently large, $\sat_{K_r}(n,F) = 0$ for all $r\geq t- \alpha + d$. 
    \begin{proof}
        Let $u$ be as defined above. Suppose $G$ is $F$-saturated with $u$ vertices of degree $n-1$. Pick $u$ such vertices and remove them one by one. Setting $\mathcal{F} = \{F\}$, we have that $G$ is $\mathcal{F}$-saturated. After repeated application of Lemma 6.2, we obtain a graph $G'$ that is $\mathcal{F'}$-saturated where $S_d \in \mathcal{F'}$. Thus the maximum degree of $G'$ is $d-1$ and $\omega(G') \leq d$ where $\omega$ is the clique number of $G$. Hence the largest clique in $G$ has size at most $t-\alpha +d-1$, and the proof is complete. 
    \end{proof}
\end{proposition}

Based on the results above for stars and paths, along with the fact that among all trees $T$, the star $S_t$ has the largest saturation number \cite{KT}, one may suspect that for any tree $T$ there exist triangle-free graphs on $n$ vertices for $n$ sufficiently large that are $T$-saturated. We note that this is true for any tree $T$ on $t\leq 6$ vertices. 
However, the following example $T^*$ shows that this does not hold for all trees.
\begin{figure}[H]
    \centering
    \begin{tikzpicture}[scale =0.7]
        \draw[fill=black] (0,0) circle (3pt);
        \draw[fill=black] (1,0) circle (3pt);
        \draw[fill=black] (2,0) circle (3pt);
        \draw[fill=black] (-2,0) circle (3pt);
        \draw[fill=black] (-1,-0) circle (3pt);
        \draw[fill=black] (0,-2) circle (3pt);
        \draw[fill=black] (0,-1) circle (3pt);
        \draw[thick] (-2,0) -- (-1,0) -- (0,0) -- (1,0) -- (2,0);
        \draw[thick] (0,0) -- (0,-1) -- (0,-2);
    \end{tikzpicture}
    \caption{$T^*$}
\end{figure}

\begin{figure}[H]
            \centering
            \begin{tikzpicture}[scale =0.6]
                \draw[fill=black] (1,-1) circle (3pt);
                \draw[fill=black] (2,-0.5) circle (3pt);
                \draw[fill=black] (2,0) circle (3pt);
                \draw[fill=black] (2,0.5) circle (3pt);
                \draw[fill=black] (1,0) circle (3pt);
                \draw[fill=black] (1,1) circle (3pt);
                \draw[fill=black] (0,2) circle (3pt);
                \draw[fill=black] (-1,1) circle (3pt);
                \draw[fill=black] (-1,0) circle (3pt);
                \draw[fill=black] (-1,-1) circle (3pt);
                \draw[fill=black] (0,-2) circle (3pt);
                \draw[fill=black] (0,-3) circle (3pt);
                \draw[fill=black] (0.5,3) circle (3pt);
                \draw[fill=black] (-0.5,3) circle (3pt);
                \draw[fill=black] (2,1.5) circle (3pt);
                \draw[fill=black] (-2,1.5) circle (3pt);
                \draw[fill=black] (-2,0) circle (3pt);
                \draw[fill=black] (2,-1.5) circle (3pt);
                \draw[fill=black] (-2,-1.5) circle (3pt);
                \draw[fill=black] (-2,-1) circle (3pt);
                \draw[thick] (1,-1) -- (1,0) -- (1,1) -- (0,2) -- (-1,1) -- (-1,0) -- (-1,-1) -- (0,-2) -- (1,-1);
                \draw[thick] (2,0) -- (1,0) -- (2,-0.5);
                \draw[thick] (1,0) -- (2,0.5);
                \draw[thick] (0,-3) -- (0,-2);
                \draw[thick] (0.5,3) -- (0,2) -- (-0.5,3);
                \draw[thick] (1,1) -- (2,1.5);
                \draw[thick] (-1,1) -- (-2,1.5);
                \draw[thick] (-1,0) -- (-2,0);
                \draw[thick] (-2,-1) -- (-1,-1) -- (-2,-1.5);
                \draw[thick] (1,-1) -- (2,-1.5);
                \draw[dashed] (2,-1.5) -- (0,-2);
        \end{tikzpicture}
        \caption{Example of a cycle with pendants. The dotted edge does not induce $T^*$ when added.}
        \label{pendant}
    \end{figure}
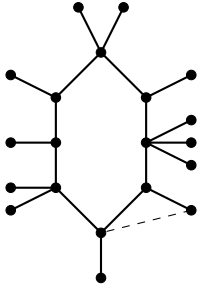

\begin{proposition}
    There does not exist a $K_3$-free graph that is $T^*$-saturated.
    \begin{proof}
        Suppose $G$ is a $T^*$-saturated graph that is $K_3$-free; clearly $n=n(G)\geq 7$. Since $G$ is $K_3$-free, the neighborhood of any vertex must be an independent set. If $G$ has a vertex of degree $2$, then adding the missing edge between its neighbors must create a copy of $T^*$. In particular, this copy must use the added edge. There are two cases. Either the added edge is incident to the degree $3$ vertex in $T^*$, or it is incident to a degree $1$ vertex in the created copy of $T^*$. In either case, we can replace the added edge with one of the original edges, contradicting the fact that $G$ is $T^*$-free. Thus $G$ has no vertex of degree $2$.  
        
        Now, the maximum degree of $G$ must be at least $3$. Otherwise our graph is a matching along with some isolated vertices, but such a graph is not $T^*$-saturated. We will only focus on a component containing a vertex of degree at least $3$ as any edge added within that component must create $T^*$. Let $x$ be a vertex of degree at least $3$. If $x$ has three neighbors $a,b,c$ of degree at least $3$, then their neighborhoods must be precisely $x$ and two other common vertices. Otherwise $G$ already contains $T^*$. There must be another vertex $u$ in this component of $G$ because it has missing edges and not enough vertices to create a copy of $T^*$. A quick check shows that no matter which vertex we join $u$ to, we must already have a copy of $T^*$, contradicting the condition that $G$ is $T^*$-saturated. 
        
        We also note that if $x$ has a unique neighbor of degree at least $3$, then adding an edge between two of its neighbors of degree $1$ will not create a copy of $T^*$. Therefore every vertex of degree at least $3$ is adjacent to at least one vertex of degree $1$ and exactly two vertices of degree at least $3$. Therefore $G$ is isomorphic to a cycle whose vertices each have at least one pendant. However, these graphs are not $T^*$-saturated. This can be seen by adding an edge between a vertex on the cycle and a pendant of a neighboring vertex as in Figure \ref{pendant}. This contradicts the only remaining case. Therefore no such graph $G$ exists.
    \end{proof}
\end{proposition}

\section{Concluding remarks}
It is still unknown what the value of $\sat_{S_2}(n,K_3)$ is. Ergemlidze, Methuku, Tait, and Timmons show in \cite{EMTT} that the split graph is not always optimal for this problem and remark that there is a connection to Moore graphs with diameter $2$ and girth $5$. It would be interesting to know what the exact value of the generalized saturation number is in this specific case. 

It would also be interesting to more closely study $\sat_{S_r}(n,S_t)$ when $t$ is even and when $t+1 \leq n \leq 2t-1$. There is less control over the structure of $S_t$-saturated graphs in this setting and more care is required. 
We are also interested in refining our bounds on $n$ for the existence of $S_t$-saturated graphs that are $K_{r+1}$-free.

\end{document}